\documentstyle{amsart}
\def\dj{d\kern-.30em\raise1.25ex\vbox{\hrule width .3em height .03em}}
\def\Dj{D\kern-.70em\raise0.75ex\vbox{\hrule width .3em height .03em}
\kern.03em}
\newsymbol\restr 1316
\newcommand{\cl}{\mbox{\family{euf}\shape{n}\selectfont cl}}

\newcommand{\id}{\mbox{\shape{n}\selectfont id}}
\newcommand{\Sum}{\displaystyle{\sum}}
\newcommand{\im}{\operatorname{im}}
\newcommand{\ita}{\iota^F}

\newtheorem{lem}{Lemma}
\newtheorem{pro}[lem]{Proposition}
\newtheorem{thm}[lem]{Theorem}
\theoremstyle{definition}
\newtheorem{defn}{Definition}
\def\bC{\Bbb{C}}
\begin{document}
\title[Braided Clifford Algebras]{Clifford Algebras and Spinors \\
for Arbitrary Braids}
\author{Mi\'co {\Dj}ur{\Dj}evi\'c \& Zbigniew Oziewicz}
\address{Instituto de Matematicas, UNAM,
Area de la Investigacion
Cientifica, Circuito Exterior, Ciudad Universitaria, M\'exico DF,
CP 04510, M\'EXICO\newline
\indent Centro de Investigaciones Teoricas, Facultad de Estudios
Superiores Cuautitlan, UNAM,
Cuautitlan Izcalli, CP 54700, M\'EXICO}
\maketitle
\begin{abstract}
General braided counterparts of classical Clifford algebras are
introduced and investigated. Braided Clifford algebras are defined as
Chevalley-K\"ahler deformations of the corresponding braided exterior
algebras. Analogs of the spinor representations are studied,
generalizing classical Cartan's approach. It is shown that, under
certain assumptions concerning the braiding, the spinor representation
is faithful and irreducible, as in the classical theory.
\end{abstract}
\section{Introduction}
\renewcommand{\thepage}{}
In this study classical theory of Clifford algebras and spinors will be
incorporated into a general braided framework. The main idea that will
be followed is the Chevalley-K\"ahler interpretation of Clifford
algebras, as deformations of exterior algebras. General braided
Clifford algebras will be introduced by constructing a new product in
the braided exterior algebras spaces.

The paper is organized as follows. In the
next section
the main properties of braided exterior algebras \cite{w} are collected.
In Section~3
we introduce and analyze a braided inner product.
In particular, we deal with both exterior algebras
over the space and its dual. These algebras are mutually dual in a natural
manner. It turns out that
the inner product with elements of the dual exterior algebra
can be viewed as the transposed operator
of the corresponding exterior multiplication, as in the classical theory.

Section 4 is devoted to the construction and general analysis of braided
Clifford algebras. Our construction conceptually follows classical
Chevalley-K\"ahler approach. We introduce a new product in the
exterior algebra space.
This product is expressible in terms of the exterior
product, and various ``relative'' contractions, which are constructed from
the corresponding scalar product in the vector space and the braided
inner product.
In such a way the Clifford algebra becomes a
deformation of the exterior algebra. We also introduce an analog
of the
Crumeyrolle map \cite{cru} in the tensor algebra, connecting Clifford
and exterior ideals.
This can be viewed as another way to define
braided Clifford algebras.

In Section 5 we study counterparts of algebraic spinors, conceptually
following Cartan's geometrical approach,
and varying and generalizing a construction given in \cite{bau}.
Spinors are defined as elements of
braided exterior algebras over certain isotropic subspaces of the
initial vector space. By construction the spinor space
is a left Clifford module. We prove that, under certain assumptions
concerning the braid, the spinor representation is irreducible and faithful,
as in the classical theory.
Finally, in Section 6 some concluding remarks are made.

\section{On Braided Exterior Algebras}
\renewcommand{\thepage}{\arabic{page}}
In this section the main properties of exterior algebras \cite{w}
associated to
an arbitrary braid operator are collected. Let $W$ be a (complex)
finite-dimensional vector space, and let $\sigma\colon W\otimes W\rightarrow
W\otimes W$ be a bijective map satisfying the braid equation
\begin{equation}\label{braid}
(\sigma\otimes\id)(\id\otimes\sigma)(\sigma\otimes\id)
=(\id\otimes\sigma)(\sigma\otimes\id)(\id\otimes\sigma).
\end{equation}

Let $A\colon W^\otimes\rightarrow W^\otimes$ be
the corresponding total antisymetrizer map. Its components $A_n\colon
W^{\otimes n}\rightarrow W^{\otimes n}$ are given by
$$ A_n=\sum_{\pi\in S_n} (-1)^\pi\sigma_\pi$$
where $\sigma_\pi\colon W^{\otimes n}\rightarrow W^{\otimes n}$ are
maps obtained by replacing transpositions
figuring in a minimal decomposition of  $\pi$
by the corresponding $\sigma$-twists. The following identities hold
\begin{align}
A_{n+k}&=\bigl(A_n\otimes A_k\bigr)A_{nk}\label{dec1}\\
A_{n+k}&=B_{nk}\bigl(A_n\otimes A_k\bigr)\label{dec2},
\end{align}
where
$$A_{nk}=\sum_{\pi\in S_{nk}}(-1)^\pi\sigma_{\pi^{-1}},\qquad
B_{nk}=\sum_{\pi\in S_{nk}}(-1)^\pi\sigma_\pi$$
and $S_{nk}\subseteq S_{n+k}$ is the set of permutations preserving the
order of sets $\{1,\ldots,n\}$ and $\{n+1,\ldots,n+k\}$.

By definition \cite{w},
the corresponding braided {\it exterior algebra}
$W^\wedge(\sigma)=W^\wedge$ is the
factoralgebra of the tensor algebra $W^\otimes$ relative
to the ideal $\ker(A)$.

The algebra $W^\wedge$ can be naturally realized as a subspace
$\im(A)$
in $W^\otimes$. This realization is given by
\begin{equation}\label{real}
\Bigl[\psi+\ker(A)\Bigr]\, \leftrightarrow\, A(\psi).
\end{equation}
In terms of the above identification the exterior product
of $\psi\in W^{\wedge n}$ and $\varphi\in W^{\wedge k}$ is given by
\begin{equation}\label{mult}
\psi\wedge\varphi=B_{nk}(\psi\otimes\varphi).
\end{equation}

\section{Braided Inner Product}
The spaces $W^{*\otimes n}$ and $W^{\otimes n}$
are mutually dual, in a natural manner. The duality between them is given
by
\begin{equation}\label{dual}
\langle f_1\otimes\dots\otimes f_n,\psi_n\otimes\dots\otimes\psi_1\rangle=
f_1(\psi_1)\dots f_n(\psi_n).
\end{equation}
This pairing trivialy extends to the whole tensor algebras.
We will assume that the dual space is equiped with the transposed
braiding
$\sigma^*\colon W^*\otimes W^*\rightarrow W^*\otimes W^*$. In what
follows operators associated to $\sigma^*$ will be endowed with the star.

Maps $\{A_n,A_n^*\}$, as well as $\{A_{nk}, B_{kn}^*\}$ and $\{A_{nk}^*,
B_{kn}\}$ are mutually transposed.
Furthermore, it is possible to define a natural
pairing $(\,)^\wedge$
between
exterior algebras $W^{*\wedge}$ and $W^{\wedge}$. Explicitly,
$$ (\phi,\vartheta)^\wedge=\langle\phi,\widetilde{\vartheta}\rangle,$$
where $W^{*\wedge}$ is embedded in $W^{*\otimes}$, and
$[\widetilde{\vartheta}]^\wedge=\vartheta$.

For each $f\in W^*$ and $\xi\in W^{\wedge n}$, let
$f\sqcup\xi\in W^{\wedge n-1}$ be an element given by
\begin{equation}\label{contr}
f\sqcup\xi=(f\otimes\id^{n-1})(\xi).
\end{equation}
In the above formula, it is assumed that $W^\wedge$ is
embedded in $W^\otimes$. The fact
that $f\sqcup\xi$ belongs to $W^\wedge$ follows from the
decomposition \eqref{dec1}.

In such a way we have constructed a map $\sqcup\colon W^*\otimes
W^\wedge\rightarrow W^\wedge$ (a counterpart
of the standard inner product with 1-forms). For each $f\in W^*$ we will
denote by $\sqcup_f\colon W^\wedge\rightarrow W^\wedge$ the corresponding
contraction map.

In what follows it will be assumed
that $\sigma$ is naturally extended to a braiding on
$W^\wedge\otimes W^\wedge$,
by requiring
\begin{align}
\sigma(m\otimes\id)&=(\id\otimes
m)(\sigma\otimes\id)(\id\otimes\sigma)\label{s-m1}\\
\sigma(\id\otimes m)&=(m\otimes\id)(\id\otimes \sigma)(\sigma\otimes
\id),\end{align}
where $m\colon W^\wedge\otimes W^\wedge\rightarrow W^\wedge$ is the
product map.

\begin{lem} The following identity holds
$$ \sqcup_f(\xi\eta)=\sqcup_f(\xi)\eta+
(-1)^{\partial\xi}
m\sigma^{-1}(\sqcup_f\otimes\id)\sigma(\xi\otimes\eta).$$
This is a counterpart of the standard graded Leibniz rule.
\end{lem}
\begin{pf}
The statement follows from the definition of $\sqcup$,
decomposition \eqref{dec1},
and correspondence \eqref{real}.
\end{pf}

The introduced operator $\sqcup$ can be trivially extended to the map
of the form
$\sqcup\colon W^{*\otimes}\otimes W^{\wedge}\rightarrow W^{\wedge}$
such that
$$\sqcup_{f\otimes g}=\sqcup_f\sqcup_g$$
for each $f,g\in W^{*\otimes}$.
\begin{lem} If $f\in\ker(A^*)$ then $\sqcup_f=0$.
\end{lem}
\begin{pf} The statement follows
from decomposition \eqref{dec1}, and the
definiton of $\sqcup$.
\end{pf}

Hence, we can pass from $W^{*\otimes}$ to $W^{*\wedge}$ in the first argument
of $\sqcup$. In such a way we obtain a map of the form
$\sqcup
\colon W^{*\wedge}\otimes W^\wedge\rightarrow W^\wedge$ (we use the same
symbol for different maps, because the domain is clear from
the context).

\begin{defn} The constructed map is called {\it the braided inner
product}.
\end{defn}

The inner product map can be equivalently described
as the transposed exterior
multiplication.
For each $f\in W^{*\wedge}$, let
$\wedge_f\colon W^{*\wedge}\rightarrow
W^{*\wedge}$ be a linear map given by $\wedge_f(\phi)=\phi\wedge f$.

\begin{pro} The following identity holds
$$ (\phi, \sqcup_f(\xi))^\wedge=(\wedge_f(\phi),\xi)^\wedge. $$
In other words, $\sqcup_f$ and $\wedge_f$ are
mutually transposed.
\end{pro}
\begin{pf}
If $\phi\in W^{*\wedge k}$ and $\xi\in W^{\wedge k+r}$ then
\begin{multline*}
(\wedge_f(\phi),\xi)^\wedge=(\phi\wedge f,\xi)^\wedge=\langle
B_{kr}^*(\phi\otimes f),\widetilde{\xi}\rangle\\
=\langle\phi\otimes f,A_{rk}\widetilde{\xi}\rangle
=(\phi, f\sqcup\xi)^\wedge=(\phi,\sqcup_f(\xi))^\wedge
\end{multline*}
for each $f\in W^{*\wedge r}$.
\end{pf}

\section{Braided Clifford Algebras}

Let us assume that $W$ is endowed with a scalar product $F$, understood
as a linear map $F\colon W\otimes W\rightarrow \bC$.
Let $\ell_F\colon W\rightarrow W^*$ be an associated correlation,
given by
$$ \bigl[\ell_F(x)\bigr](y)=F(x\otimes y).$$

As first, we are going to introduce various types of
contraction operators in $W^\wedge$
which will play a fundamental role in constructing the corresponding
Clifford algebra.

In what follows it will be
assumed that $F$ and $\sigma$ are mutually related such that the
following ``funny functoriality'' holds
\begin{equation}\label{ffun}
(F\otimes \id)(\id\otimes\sigma)=(\id\otimes F)(\sigma\otimes\id).
\end{equation}
This implies
\begin{equation}
(\ell_F\otimes\ell_F)\sigma=\sigma^*(\ell_F\otimes\ell_F).
\end{equation}
In particular
\begin{equation}
\ell_F^\otimes\bigl(\ker(A)\bigr)\subseteq\ker(A^*),
\end{equation}
where $\ell_F^\otimes$ is the unital multiplicative extension of
$\ell_F$. Factorizing $\ell_F^\otimes$ through ideals
$\ker\{A,A^*\}$
we obtain a homomorphism $\ell_F^\wedge\colon
W^\wedge\rightarrow W^{*\wedge}$.

Let $\ita\colon W^\wedge\otimes W^\wedge\rightarrow W^\wedge$ be
a contraction map given by
\begin{equation}
\ita=\sqcup(\ell_F^\wedge\otimes\id).
\end{equation}

By construction, $\ita$ is multliplicative
on the first factor and satisfies the following braided
variant of the Leibniz rule
$$\ita_x(\vartheta\eta)=\ita_x(\vartheta)\eta+(-1)^{\partial\vartheta}
\sum_k\vartheta_k\ita_{x_k}(\eta)$$
where $x\in W$ and
$\Sum_k\vartheta_k\otimes x_k=\sigma(x\otimes\vartheta)$.

Let us define relative contraction operators
$\langle\,\rangle_k\colon W^\wedge\times W^\wedge\rightarrow W^\wedge$ as
follows
$$\langle\zeta,\xi\rangle_k=
\sum_j\psi_j\wedge\bigl(\ita_{\varphi_j}\xi\bigr)$$
Here, it is assumed that $\zeta\in W^{\wedge n}$ and
$\Sum_j\psi_j\wedge\varphi_j=[A_{n-kk}\widetilde{\zeta}]^\wedge$,
where $\varphi_j\in W^{\wedge k}$, $\psi_j\in W^{\wedge n-k}$
and $\widetilde{\zeta}\in W^{\otimes n}$ satisfies
$[\widetilde{\zeta}]^\wedge=\zeta$,.
Consistency of this definition follows from \eqref{dec1}. If $n<k$ we
define $\langle\,\rangle_k=0$.

Now, we define a {\it new product} on $W^\wedge$, in the spirit
of the classical construction. This product is
introduced by the following expression
$$\psi\circ\varphi=\psi\wedge\varphi+\sum_{k\geq 1}
\langle\psi,\varphi\rangle_k.$$
In particular,
$$x\circ\psi=x\wedge\psi+\ita_x(\psi)$$
for $x\in W$.
This generalizes classical Chevalley's formula.
\begin{thm} Endowed with $\circ$, the space $W^\wedge$ becomes
a unital associative algebra, with the unity $1\in W^\wedge$.
\end{thm}
\begin{pf}
The standard diagramatic computations, using braid
diagrams and functoriality property \eqref{ffun}.
\end{pf}

\begin{defn} The algebra $\cl(W,\sigma,F)=(W^\wedge,\circ)$ is called {\it
the braided Clifford algebra} (associated to $\{W,\sigma,F\}$).
\end{defn}

The constructed algebra can be understood as a deformation of
the exterior algebra $W^\wedge$. The graded algebra associated
to the filtered algebra $\cl(W,\sigma,F)$ naturally coincides with $W^\wedge$.

The algebra $\cl(W,\sigma,F)$ can be viewed as a factoralgebra
$\cl(W,\sigma,F)=W^\otimes/J_F$,
where $J_F$ is the kernel of the canonical epimorphism $j_F\colon W^\otimes
\rightarrow \cl(W,\sigma,F)$ (extending the identity map on $W$). Now,
this ideal will be described in an independent
way, using a generalization of the construction presented
by Crumeyrolle in \cite{cru}.

As first, a Clifford product in the tensor algebra $W^\otimes$
will be introduced.
Let us consider a linear map $\lambda_F\colon W^\otimes\rightarrow
W^\otimes$ defined by
$$\lambda_F(1)=1\qquad\lambda_F(x\otimes\vartheta)=x\otimes
\lambda_F(\vartheta)+\ita_x\lambda_F(\vartheta)$$
where $x\in W$ and $\vartheta\in W^\otimes$.
In the above formula, $\ita_x$ is considered as a braided antiderivation
on $W^\otimes$.
The map $\lambda_F$ is bijective. Let
$\circ$ be a new product in $W^\otimes$, given by
$$\vartheta\circ\eta=\lambda_F\bigl(\lambda_F^{-1}(\vartheta)
\otimes\lambda_F^{-1}(\eta)\bigr).$$

By construction the space $\ker(A)$ is a left ideal in $W^\otimes$,
relative to this new product. Condition \eqref{ffun}
ensures that $\ker(A)$ is
also a right $\circ$-ideal.
\begin{thm} We have
$$(W^\otimes,\circ)/\ker(A)=\cl(W,\sigma,F).$$
\end{thm}
\begin{pf}
It is sufficient to observe that
\eqref{ffun} implies that the product $\circ$ on $W^\otimes$ is given by
essentially the same formula as for $W^\wedge$, the only difference is that
$\wedge$ should be replaced by $\otimes$, and contractions are acting
on $W^\otimes$.
\end{pf}

In other words, the factorization map $[\,]^\wedge\colon
W^\otimes\rightarrow W^\wedge$ is also a homomorphism
of corresponding deformed algebras. The map $\lambda_F$ is a
braided counterpart of the maps introduced in \cite{cru} and
\cite{ozirot}.
\begin{lem} We have
\begin{equation}\label{J-Ext}
\ker(A)=\lambda_F(J_F).
\end{equation}
\end{lem}
\begin{pf}
The statement follows from the equality
$[\,]^\wedge\lambda_F=j_F$.
\end{pf}
\section{The Spinor Representation}

This section is devoted to a braided generalization of classical
Cartan theory of spinors \cite{car}. We follow ideas of \cite{bau}.
Let us assume that the space $W$ is splitted
into a direct sum
$$W=W_1\oplus W_2$$
where $W_1$, $W_2$ are $F$-isotropic subspaces. Furthermore, let
us assume that this decomposition is compatible with the braiding $\sigma$
in the following way
\begin{gather}
\sigma(W_i\otimes W_j)=W_j\otimes W_i\label{S-dec1}\\
\sigma^2\Big\vert\Bigl\{(W_1\otimes W_2)\oplus(W_2\otimes
W_1)\Bigr\}=\id.\label{S-dec2}
\end{gather}

Finally, it will be assumed that $F\restr \bigl(W_1\otimes W_2\bigr)=0$
and that $F\restr \bigl(W_2\otimes W_1\bigr)$ is {\it nondegenerate}. In this
case,
$W_2=W_1^*$, in a natural manner. The duality is given by
$F(f\otimes x)=f(x)$,
for $f\in W_2$ and $x\in W_1$.

Exterior algebras $W_1^\wedge$ and $W_2^\wedge$  are understandable
as subalgebras of $\cl(W,\sigma,F)$, in a natural manner.
\begin{lem}\label{lem:7} The map $\mu\colon W_1^\wedge\otimes W_2^\wedge
\rightarrow \cl(W,\sigma,F)$ defined by
\begin{equation}\label{car}
\mu(u\otimes v)=u\circ v
\end{equation}
is bijective.
\end{lem}
\begin{pf}
If $u\in W_1$ and $v\in W_2$ then
$$vu+\sum_k u_k v_k-F(v\otimes u)1=0$$
where $\Sum_ku_k\otimes v_k=\sigma(v\otimes u)$. This implies that $\mu$
is surjective.

Let $\mu^*\colon W_1^\wedge\otimes W_2^\wedge
\rightarrow W^\wedge$ be the grade-preserving component of $\mu$. This map is
explicitly given by $\mu^*(u\otimes v)=u\wedge v$. We prove that
$\mu^*$ is injective. We have
$$\mu^*(u\otimes v)=B_{kl}(u\otimes v)$$
for $u\in W_1^{\wedge k}$ and $v\in W_2^{\wedge l}$. If $\psi\in
W_1^{\wedge k}\otimes W_2^{\wedge l}$ then $p_{kl}\mu^*(\psi)=\psi$,
where $p_{kl}\colon W^{\otimes k+l}\rightarrow W_1^{\otimes k}\otimes
W_2^{\otimes l}$ is the projection map. Hence $\mu^*$ is injective.
\end{pf}

Let us consider the space ${\cal{K}}=W_2^\wedge$, and let
$\kappa
\colon {\cal{K}}\rightarrow\bC$ be a natural character, specified by
$\kappa(1)=1$ and $\kappa(W_2)=\{0\}$. This gives a left ${\cal{K}}$-module
structure on the number field $\bC$. On the other hand, $\cl(W,\sigma,F)$ is a
right ${\cal{K}}$-module,
in a natural manner. Let ${\cal{S}}$ be a left $\cl(W,\sigma,F)$-module, given
by
$$ {\cal{S}}=\cl(W,\sigma,F)\otimes_{\cal{K}}\bC.$$
\begin{defn} The constructed module is called {\it the spinor module}
for $\cl(W,\sigma,F)$ associated to $\bigl\{W_1,W_2\bigr\}$.
\end{defn}
According to Lemma~\ref{lem:7}, the space $\cal{S}$ is  naturally identificable
with the exterior algebra $W_1^\wedge$. In terms of this identification, we
have
$$ x\xi=x_1\wedge\xi+x_2\sqcup\xi,$$
for each $x\in W$, where $x=x_1+x_2$  and $x_i\in W_i$. In other words,
a complete analogy with the classical Cartan formalism holds. The
formula \eqref{car} corresponds to the classical Cartan map.

\begin{thm} The algebra $\cl(W,\sigma,F)$ acts on $\cal{S}$
faithfully and irreducibly.
\end{thm}
\begin{pf}
We first prove that each vector
$\psi\in{\cal{S}}\setminus\{0\}$
is cyclic (which implies that the module is simple).
Obviously, the unit element $1_{\cal{S}}\in\cal{S}$ is
cyclic, by construction. The duality between spaces $W_1$ and $W_2$ extends
to the duality between exterior algebras $W_1^\wedge$ and $W_2^\wedge$, as
explained in Section~3. In terms of this duality, the contraction between
elements of the same degree is just the pairing map. It follows that
there exists $\varphi\in{\cal{K}}$ such that $\varphi\psi=1_{\cal{S}}$. Hence,
$\psi$ is cyclic.

Let us consider an element $x\in\cl(W,\sigma,F)\setminus\{0\}$. We have
$$\mu^{-1}(x)=\sum_k u_k\otimes v_k+\psi$$
where $\Sum_k u_k\otimes v_k\neq 0$ is the component consisting of
summands having the minimal second degree $n$. We can assume that $v_k$
are linearly independent vectors.

There exist spinors $\xi_j\in\cal{S}$ satisfying
$\psi\xi_j=0$ and $v_k\xi_j=\delta_{kj}$.
This gives $x\xi_j=u_j$ which implies that
the representation is faithful.
\end{pf}

The module $\cal{S}$ is completely characterized by the existence of
a cyclic vector (the unit element $1_{\cal{S}}$), killed by the space
$W_2$.

In other words let $\cal{V}$ be an arbitrary
(left) $\cl(W,\sigma,F)$-module, possesing a
vector $v$ satisfying $\bigl\{W_2\bigr\}v=\{0\}$. Then there exists
the unique module map $\varrho\colon\cal{S}\rightarrow\cal{V}$
satisfying $\varrho\bigl(1_{\cal{S}}\bigr)=v$. The map $\varrho$ is
injective (because of the simplicity of $\cal{S}$). In particular, if
$v$ is cyclic then $\varrho$ is a module isomorphism.

\section{Concluding Remarks}
If the braid operator $\sigma$ is such that $\ker(A)$ is quadratic, then
the ideal $J_F$ is generated by elements of the form
\begin{equation}\label{quad}
Q=\psi-F(\psi)1\otimes 1
\end{equation}
where $\psi\in W^{\otimes 2}$ is $\sigma$-invariant. This covers
Clifford algebras based on Hecke braidings \cite{ozi,ozirot}
(including classical Weyl algebras, in the trivial way).

Quantum Clifford algebras and spinors introduced and analyzed in
\cite{bau} can be included in the theory presented here. Clifford
algebras of \cite{bau} are based on Hecke braidings $\tau\colon V\otimes
V\rightarrow V\otimes V$ (where $V$ is a finite-dimensional vector
space) admitting extensions to all possible braidings
between $V$ and $V^*$, so that the contraction map is functorial, in the
standard sense. Then $W=V\oplus V^*$, the corresponding scalar product
$F$ and the braiding $\sigma$ are expressible in terms of the
extended braiding $\tau$ and the contraction map.

In the classical theory spinors can be equivalently
viewed as elements of the left $\cl(W,\sigma,F)$-ideal, generated by a volume
element of $W_2$. A similar description is possible in the braided
context, if the external algebra $W_2^\wedge$ admits ``volume
elements''. Namely let us assume that $\omega\in W_2^\wedge$ is such
that
$$\bigl\{W_2\bigr\}\omega=\{0\}.$$
Then the left $\cl(W,\sigma,F)$-ideal $\cal{I}_\omega$
generated by $\omega$ is canonically isomorphic to $\cal{S}$, as a left
$\cl(W,\sigma,F)$-module.

The construction of the map $\lambda_F$ works for an arbitrary $F$ (and
in particular, it is independent of functoriality-type
assumptions \eqref{ffun}). For a possibility to define braided Clifford
algebras as deformations of braided exterior algebras, it is sufficient
to assume that $\ker(A)$ is also a right-ideal in $(W^\otimes,\circ)$.
This assumption is weaker then \eqref{ffun}. However, if \eqref{ffun}
does not hold, then the symmetry between left and right is broken.

\end{document}